\documentclass[11pt]{article}
\usepackage[english]{babel}

\usepackage{amsmath}
\usepackage{amssymb}

\newcommand{\R}{   {\ifmmode{{\mathbb R}}\else{$\mathbb R$}\fi}}
\newcommand{\N}{   {\ifmmode{{\mathbb N}}\else{$\mathbb N$}\fi}}

\newcommand{\Li}{\mathcal{L}}
\newcommand{\Lo}{\mathcal{L}_0}

\newcommand{\cc}{\mathcal{C}}
\newcommand{\cci}{\mathcal{C}^{\infty}}
\newcommand{\X}{\mathcal{X}}
\newcommand{\U}{\mathcal{U}}
\newcommand{\BU}{\mathbf{U}}
\newcommand{\A}{\mathcal{A}}

\newcommand{\rank}{\mbox{rank}\ }
\newcommand{\Int}{\mbox{Int}}

\begin{document}

\title{Finite Time and Exact Time Controllability on Compact Manifolds}
\author{Philippe JOUAN\footnote{LMRS, CNRS UMR 6085, Universit\'e
    de Rouen, avenue de l'universit\'e BP 12, 76801
    Saint-Etienne-du-Rouvray France. E-mail: Philippe.Jouan@univ-rouen.fr}}

\date{\today}

\maketitle

\begin{abstract}
It is first shown that a smooth controllable system on a compact manifold is finite time controllable. The technique of proof is close to the one of Sussmann's orbit theorem, and no rank condition is required.

This technique is also used to give a new and elementary proof of the equivalence between controllability for essentially bounded inputs and for piecewise constant ones.

Two sufficient conditions for controllability at exact time on a compact manifold are then stated.

Some applications, in particular to linear systems on Lie groups, are provided.

\end{abstract}

\section{Introduction.}

The aim of this paper is to study the finite time and exact time controllability properties of systems on compact manifolds described by a controlled equation of the form
$$
(\Sigma) \qquad\qquad\qquad \dot{x}=f(x,u)
$$
(see Section \ref{Preliminaries} for more precise assumptions). Throughout the paper the system will be assumed to be controllable. Indeed the purpose is not here to provide sufficient conditions for controllability, but first to prove that a controllable system is finite time controllable as soon as the state space is compact, and secondly to give sufficient conditions for exact time controllability.

Such results were proved for right invariant systems on compact Lie groups by Jurdjevic and Sussmann in \cite{JS272}. We consider here the general case where the system has no particular property (controllability excepted).

\vskip 0.2cm

The first result (Section \ref{tempsfini}, Theorem \ref{finitetime}) is:

\textit{If the state space $N$ is compact and the system controllable then it is finite time controllable, that is there exists $T>0$ such that any point $x$ can be steered to any point $y$ in time less than or equal to $T$.}

Notice that no rank condition is required. In fact if this last is satisfied, even only at one point, then the proof is rather straightforward. If not the main part of the proof is contained in the one of Lemma \ref{mainlemma}, which asserts:

\textit{If the system is controllable, then for any $x\in N$ there exists $t>0$ such that the interior of $\A(x,\leq t)$ is not empty.}

This lemma is part of a more general result known as ``accessibility implies normal accessibility'' which can be found in \cite{Grasse84} (see also \cite{Sussmann76}). Lemma \ref{mainlemma} is less general but its proof is different and simpler than the original ones. In \cite{Grasse84} and \cite{Sussmann76}, it is roughly speaking proved that if a point does not have the normal accessibility property then its attainable set is included in a subset of first category. On the other hand the proof of Lemma \ref{mainlemma} is analogous to the one of the Orbit's Theorem (see \cite{Sussmann73}), with the constraint that the trajectories cannot be followed backward.

This proof does not involve compactness of the state space $N$, and can therefore be applied to general manifolds (see Theorem \ref{noncompact}).
\vskip 0.1cm
Moreover it allows as well to give a new and elementary proof of the following known result: {\it a system is controllable for locally essentially bounded inputs if and only if it is controllable for piecewise constant inputs} (Theorem \ref{piecewise} in Section \ref{equivalence} and \cite{GS90} for the original statement).
\vskip 0.2cm

Section \ref{tempsexact} is then devoted to exact time controllability on compact manifolds, and two sufficient conditions are provided. More accurately:

{\it The manifold $N$ is assumed to be compact and the system to be controllable. If one of the following conditions holds

(i) The zero-time ideal has full rank at one point (Theorem \ref{exactime}),

(ii) There exist $x\in N$ and $u\in \U$ such that $f(x,u)=0$ (Theorem \ref{singularity}),

then it is exact time controllable. If the system is Lie-determined then the first condition is necessary.}

\vskip 0.2cm

Section \ref{Applications} contains some applications, in particular to linear systems on compact Lie groups. They are systems whose drift vector field is the infinitesimal generator of a one parameter group of automorphisms, and whose controlled vector fields are invariant (see \cite{JouanAX08}). They are controllable if and only if they satisfy the rank condition, and we show that they are under that condition exact time controllable.


\section{Preliminaries}\label{Preliminaries}

We deal with the control system
$$
(\Sigma) \qquad\qquad\qquad \dot{x}=f(x,u)
$$
where $x$ belongs to a $\cci$, connected, $n$-dimensional manifold $N$, and $u$ to a subset $\U$ of $\R^m$. The state space $N$ is not required to be compact in general. The controlled vector field $f$ is assumed to be $\cci$ with respect to $(x,u)$ (if $\U$ is not an open set, this means that $f$ can be extended to a $\cci$ controlled vector field defined on $N\times V$ where $V$ is an open neighborhood of $\U$).

\vskip 0.2cm

A measurable function $\omega$ from $[0,+\infty[$ into $\U$ belongs to $L^{\infty}_{loc}([0,+\infty[,\U)$ if for all $T>0$, $\omega(t)$ belongs to some compact $K\subset \U$ for almost every $t\in [0,T]$.

The set $\BU$ of admissible inputs is a subset of $L^{\infty}_{loc}([0,+\infty[,\U)$ which contains the piecewise constant inputs and is stable with respect to concatenation, that is if $\omega$ and $\nu$ belong to $\BU$, then the function $w$ defined by 
$$
w(t)=\left\{ \begin{array}{ll}
\omega(t) & t\in [0,T[\\
\nu(t-T) & t\in [T,+\infty[
\end{array}\right .
$$
belongs as well to $\BU$.

\vskip 0.2cm

In case where $\BU$ is the set of piecewise constant inputs, the regularity assumption of the vector field with respect to the control may be relaxed, it is enough to require $f$ to be $\cci$ w.r.t. $x$ for each $u$. In other words we can consider $\{f(.,u);\ u\in\U\}$ as a family a vector fields.

\vskip 0.2cm

Given a point $x\in N$ and an admissible input $\omega\in\BU$ the trajectory of $(\Sigma)$ starting from $x$ is denoted by
$$
t\longmapsto \phi(t,x,\omega).
$$
It is defined on an interval $[0,b[$, where $b>0$ or $b=+\infty$.

\vskip 0.2cm

For $x\in N$ and $t\geq 0$, we denote by $\A(x,t)$ (resp. $\A(x,\leq t)$) the set of points of $N$ that can be reached from $x$ at time $t$ (resp. in time less than or equal to $t$). The system is said to be
$$
\begin{array}{l}
\mbox{controllable if } \qquad\qquad\qquad \forall x\in N \quad \A(x)=\bigcup_{t\geq 0}\A(x,t)=N\\
\mbox{finite time controllable if }\  \exists T>0 \mbox{ such that } \forall x\in N \quad \A(x,\leq T)=N\\ 
\mbox{exact time controllable if }\  \exists T>0 \mbox{ such that } \forall x\in N \quad \A(x,T)=N
\end{array}
$$

The interior of the reachable set $\A(x,t)$ will be denoted by $\Int(\A(x,t))$. Notice that according to the concatenation property of $\BU$, the following implication
$$
y\in \Int(\A(x,t)) \quad \mbox{and} \quad z\in \A(y,s)\Longrightarrow z\in \Int(\A(x,t+s))
$$
holds, and that similar definition and implication hold for reachability in time less than or equal to $t$.
\vskip 0.2cm

Let $\Li$ stand for the Lie algebra of vector fields generated by the family $\{f(.,u);\ u\in\U\}$, and $\Lo$ for the ideal of $\Li$ generated by the differences
$$
f(.,u)-f(.,v) \qquad\qquad u,v\in \U.
$$
They are related by the equality
$$
\Li=\R f(.,u)+\Lo
$$
which holds for any $u\in\U$.

The rank of $\Li$ (resp. $\Lo$) at a point $x$ is the dimension of the subspace $\Li(x)=\{X(x);\ X\in \Li\}$ (resp. $\Lo(x)=\{X(x);\ X\in \Lo\}$) of $T_xN$.

Let us recall that the system is said to be \textit{Lie-determined} if at each point $x\in N$, the rank of $\Li$ is equal to the dimension of the orbit of $\Sigma$ through $x$ (see for instance \cite{Jurdjevic97}).

\vskip 0.2cm

The forthcoming proofs make use of the time-reversed system. It will be referred to as $\Sigma^-$ and the various reachability sets denoted by $\A^-(x,t)$, $\A^-(x,\leq t)$, and so on. Note that $\A^-(x,t)$ (resp. $\A^-(x,\leq t)$) is the set of points of $N$ that can be steered to $x$ at time $t$ (resp. in time less than or equal to $t$).


\section{Finite time controllability}\label{tempsfini}

\newtheorem{finitetime}{Theorem}
\begin{finitetime}\label{finitetime}
The manifold $N$ is assumed to be compact and the system $(\Sigma)$ to be controllable. Then there exists $T>0$ such that $(\Sigma)$ is controllable in time less than or equal to $T$.
\end{finitetime}

The proof makes essentially use of the following lemma, which does not involve compactness:

\newtheorem{mainlemma}{Lemma}
\begin{mainlemma}\label{mainlemma}
If $(\Sigma)$ is controllable, then for any $x\in N$ there exists $t>0$ such that the interior of $\A(x,\leq t)$ is not empty.
\end{mainlemma}

\noindent\textit{Proof of Lemma \ref{mainlemma}}

The first step consists in the construction of an integrable distribution.
For each $\omega\in \BU$ and $t\geq 0$, the mapping
$$
x\longmapsto \phi(t,x,\omega)
$$
is a local diffeomorphism defined on a open subset (possibly empty) of $N$. The pseudo semigroup of such local diffeomorphisms will be denoted by $SG$. Let $D$ stand for the following family of $\cci$ vector fields
$$
D=\{\Phi_*f(.,u);\ \Phi\in SG \ \ u\in \U\}.
$$ 
As well as the elements of $SG$ these vector fields are defined on open subsets of $N$ that need not be equal to the entire manifold. Consider the distribution $\Delta$ spanned by $D$:
$$
\Delta(x)=Sp\{X(x);\ X\in D\}.
$$
First of all the rank of $\Delta$ is constant. Indeed for all $x,y$ belonging to $N$ the controllability assumption implies that there exists $\Phi\in SG$ such that $y=\Phi(x)$. Because of the construction of $\Delta$, and the fact that $SG$ is a (pseudo) semigroup, we have
\begin{equation}
\label{inclusion} T_x\Phi(\Delta(x))\subseteq \Delta(y)
\end{equation}
hence $\dim \Delta(x)\leq \dim \Delta(y)$. But $x$ and $y$ are arbitrary, and this yields that the rank of $\Delta$ is constant over $N$.

Since $\rank \Delta(x)=\rank \Delta(y)$, the inclusion (\ref{inclusion}) is in fact an equality, and it follows that
$$
\Delta(x)=(T_x\Phi)^{-1}(\Delta(y))=T_x\Phi^{-1}(\Delta(y))
$$
for all $x,y$ and $\Phi$ for which this makes sense. Consequently the distribution $\Delta$ is invariant under the action of the (pseudo) group $G$ of  local diffeomorphisms generated by $SG$. 

We are now in a position to show that the distribution $\Delta$ is $D$-invariant: let $X\in D$, and $(\Psi_t)_{t\in\R}$ its flow. By definition of $D$ there exists $\Phi\in SG$ and $u\in \U$ such that $X=\Phi_*f(.,u)$. Then the flow of $X$ is given by
$$
\Psi_t=\Phi\circ \phi^u_t \circ \Phi^{-1}
$$
where $(\phi^u_t)_{t\in\R}$ stands for the flow of $f(.,u)$. But $\Phi$, $\Phi^{-1}$ and $\phi^u_t$ belong to $G$, hence for all $t\in\R$, $\Psi_t$ belongs as well to $G$.

According to Theorem 4.2 of \cite{Sussmann73} the distribution $\Delta$ is integrable. But each attainable set is included in an integral manifold of this distribution, and the controllability assumption implies the fundamental consequence that the rank of $\Delta$ is everywhere full.

\vskip 0.2cm

Let $x\in N$. We are going to build a submanifold of $N$ included in $\A(x,\leq t)$, for some $t>0$, whose dimension is equal to the rank of $\Delta$. This last being full this submanifold will be an open subset of $N$ and the proof finished.

For $i=1,\dots ,n$ let $x_i\in N$, $\Phi_i\in SG$, and $u_i\in \U$ such that
\begin{enumerate}
\item
$x_1=x$, for $i=1,\dots ,n-1,$ $\Phi_i(x_{i+1})=x_i$, and $\Phi_n(x_1)=x_n$;
\item
the vectors $f(x_1,u_1)$ and $T_{x_2}\Phi_1\dots \ T_{x_i}\Phi_{i-1}.f(x_i,u_i)$ , $i=2,\dots,n$, form a basis of $T_xN$.
\end{enumerate}
Notice first that the points $x_i$ are not required to be distinct, and consequently that the local diffeomorphisms $\Phi_i$ are possibly equal to the identity mapping. Secondly such a sequence exists because $\Delta(x)=T_xN$. Indeed let us assume the construction made for $x_1,\dots,x_k$ and for the corresponding $u_i$ and $\Phi_i$. Let $V$ be the subspace of $T_xN$ spanned by the vectors $f(x_1,u_1)$ and $T_{x_2}\Phi_1\dots \ T_{x_i}\Phi_{i-1}.f(x_i,u_i)$ for $i=2,\dots ,k$, and let
$$
W=(T_{x_2}\Phi_1\dots \ T_{x_k}\Phi_{k-1})^{-1}(V).
$$
The dimension of $W$ is equal to $k$ and, if $k<n$, we can find $x_{k+1}$, $\Phi_k$ and $u_{k+1}$ such that $x_k=\Phi_k(x_{k+1})$ and $T_{x_{k+1}}\Phi_k.f(x_{k+1},u_{k+1})\notin W$. If this were not possible we would have $\Delta(x_k)\subseteq W\subsetneq T_{x_k}N$, a contradiction.

Consider now the mapping $\Theta$
$$
(s_1,\dots,s_n)\longmapsto \phi_{s_1}^{u_1}\circ \Phi_1\circ\dots\circ\phi_{s_{n-1}}^{u_{n-1}}\circ \Phi_{n-1}\circ\phi_{s_n}^{u_n}\circ \Phi_n(x)
$$
where $\phi_{s_i}^{u_i}$ stands for the flow of $f(.,u_i)$. It is defined on an open neighborhood of $0$ in $\R^n$, and its rank at $0$ is equal to $n$. Therefore it exists $\epsilon>0$ such that
$$
\Omega=\Theta(]0,\epsilon[^n)
$$
is a open subset of $N$. But for each $i=1,\dots, n$ there exist $\omega_i\in \BU$ and $t_i\geq 0$ such that $\Phi_i=\phi(t_i,.,\omega_i)$. This proves that $\Omega$ is included in $\A(x,\leq t)$ for
$$
t=n\epsilon+\sum_1^n t_i
$$
and ends the proof.

\hfill $\Box$

\vskip 0.2cm

\noindent\textit{Proof of Theorem \ref{finitetime}}

Let $x\in N$ and $t>0$ such that the interior of $\A(x,\leq t)$ is not empty, and let $y\in \Int(\A(x,\leq t))$. For any $z\in N$ there exists $s\geq 0$ such that $z\in \A(y,s)$, hence $z\in \Int(\A(x,\leq t+s))$. Therefore
\begin{equation}
\label{union}   N=\bigcup_{\tau>0} \Int(\A(x,\leq \tau))
\end{equation}
This is an increasing union of open sets, and $N$ being assumed to be compact, there exists $T_1>0$ such that $N=\Int(\A(x,\leq T_1))$, hence $N=\A(x,\leq T_1)$.

Now $\Sigma^-$ is also controllable, and by the same reasoning there exists $T_2>0$ such that $N=\A^-(x,\leq T_2)$. The system is clearly controllable in time less than or equal to $T=T_1+T_2$.

\hfill $\Box$

\vskip 0.2cm

\noindent{\bf Remark}. As soon as the system is controllable Formula (\ref{union}) is true without assuming the state space to be compact. Hence any compact subset $K$ of $N$ is covered by the increasing union of open sets $\Int(\A(x,\leq \tau))$, and, as in the proof of Theorem \ref{finitetime}, there exists a positive time $t$ such that $K$ is included in the interior of $\A(x,\leq t)$. We can therefore state:

\newtheorem{noncompact}[finitetime]{Theorem}
\begin{noncompact}\label{noncompact}
If $(\Sigma)$ is controllable then for all $x\in N$
$$
N=\bigcup_{\tau>0} \Int(\A(x,\leq \tau))
$$
and any compact subset $K$ of $N$ is included in $\Int(\A(x,\leq t))$ for some $t>0$.

\end{noncompact}


\section{Controllability for $L^{\infty}$ and for piecewise constant inputs}\label{equivalence}

In \cite{GS90} the authors proved that a system is controllable for $L^{\infty}$ inputs if and only if it is controllable for piecewise constant inputs. However Lemma \ref{mainlemma} allows to give a new, elementary, and geometric proof.

Let $\mathcal{PC}(\U)$ stand for the set of piecewise constant inputs with values in $\U$, and $\mathcal{PC}([0,t],\U)$ for the set of their restrictions to $[0,t]$. The attainable set from $x$ and for $\BU=\mathcal{PC}(\U)$ will be denoted by $\A_{pc}(x)$. Let us recall that (see for instance \cite{Sontag98}):
\begin{enumerate}
 \item $\mathcal{PC}([0,t],\U)$ is dense in $L^{\infty}([0,t],\U)$ for the $L^{\infty}$ norm.
\item For all $t>0$ the mapping
$$
(x,\omega)\in N\times L^{\infty}([0,t],\U)\longmapsto \phi(t,x,\omega)
$$
is continuous on its domain.
\end{enumerate}
From these properties it is clear that $\A_{pc}(x)$ is dense in $N$ for all $x$ as soon as $(\Sigma)$ is controllable for $\BU=L^{\infty}_{loc}([0,+\infty[,\U)$.

Let us now go back to the proof of Lemma \ref{mainlemma}. Its conclusion followed from the fact that the mapping $\Theta$
$$
(s_1,\dots,s_n)\longmapsto \phi_{s_1}^{u_1}\circ \Phi_1\circ\dots\circ\phi_{s_{n-1}}^{u_{n-1}}\circ \Phi_{n-1}\circ\phi_{s_n}^{u_n}\circ \Phi_n(x)
$$
is defined on an open neighborhood of $0$ in $\R^n$, and has full rank at this point. Recall that $\phi_{s_i}^{u_i}$ stands for the flow of the vector field $f(.,u_i)$, and that $\Phi_i$ belongs to SG. This means that for each $i=1,\dots, n$ there exist $\omega_i\in \BU$ and $t_i\geq 0$ such that $\Phi_i=\phi(t_i,.,\omega_i)$. For each $i=1,\dots, n$ pick a sequence $(\omega_i^k)_{k\geq 1}$ of piecewise constant inputs defined on $[0,t_i]$ and that converges to $\omega_i$ on this interval for the $L^{\infty}$ norm.
For $k$ large enough $\Phi_i^k:=\phi(t_i,.,\omega_i^k)$ is defined in a neighborhood of $x_{i+1}$ (of $x_1$ if $i=n$) and
$$
\Phi_i^k(x_{i+1})\longmapsto_{k \mapsto\infty} x_i.
$$
The mappings
$$
\Theta^k\ : \qquad (s_1,\dots,s_n)\longmapsto \phi_{s_1}^{u_1}\circ \Phi_1^k\circ\dots\circ\phi_{s_{n-1}}^{u_{n-1}}\circ \Phi_{n-1}^k\circ\phi_{s_n}^{u_n}\circ \Phi_n^k(x)
$$
are consequently well defined for $k$ large enough. Moreover the mapping $(x,\omega)\longmapsto\phi(t,x,\omega)$ being smooth (see \cite{Sontag98}) the differential of $\Theta$ with respect to $(s_1,\dots,s_n)$ depends continuously on the inputs $(\omega_1,\dots,\omega_n)$. Therefore the rank of $\Theta^k$ at $0$ is full for $k$ large enough.
This implies that the interior of $\A_{pc}(x)$ is not empty.

Let us choose $y$ in the interior of $\A_{pc}(x)$.  As $\A_{pc}(y)$ is dense in $N$ we conclude at once that $\A_{pc}(x)$ contains an open and dense subset of $N$. The same conclusion holds for the time reversed system: for all $x,y\in N$ the interiors of $\A_{pc}(x)$ and $\A^-_{pc}(y)$ are dense in $N$. Their intersection is consequently not empty and $y\in \A_{pc}(x)$. We have proved

\newtheorem{piecewise}[finitetime]{Theorem}
\begin{piecewise}\label{piecewise}
A $\cci$ system on a connected finite dimensional manifold is controllable for locally essentially bounded inputs if and only if it is controllable for piecewise constant inputs.
\end{piecewise}


\section{Exact time controllability}\label{tempsexact}

\newtheorem{exactime}[finitetime]{Theorem}
\begin{exactime}\label{exactime}
The manifold $N$ is assumed to be compact and the system $(\Sigma)$ to be controllable. If there exists one point $x$ where $\rank(\Lo)(x)=\dim N$, then there exists $t_0>0$ such that $(\Sigma)$ is controllable at exact time $t$ for all $t\geq t_0$.

If the system is Lie-determined, then this condition is also necessary.
\end{exactime}

As in Section \ref{tempsfini} the proof of this theorem makes use of a lemma which does not involve compactness.

\newtheorem{Thelemma}[mainlemma]{Lemma}
\begin{Thelemma}\label{Thelemma}
Let $(\Sigma)$ be a controllable system and let $x$ be a point of $N$ that verifies $\rank(\Lo)(x)=\dim N$. Then there exists $S>0$ such that
$$
\forall t\geq S \qquad \qquad x\in \ \Int(\A(x,t)).
$$
\end{Thelemma}

\noindent\textit{Proof of Lemma \ref{Thelemma}}

The assumption $\rank(\Lo)(x)=\dim N$ implies that the interiors of $\A(x,t)$ and $\A^-(x,t)$ are not empty for
$t>0$ sufficiently small (see for instance \cite{Jurdjevic97}, proof of Theorem 3, page 71). Let 
$$
x_1\in \Int(\A(x,t)) \qquad \mbox{ and } \qquad x_2\in \Int(\A^-(x,t)).
$$
Since the system is controllable, there exists a time $t'\geq 0$ such that $x_2\in \A(x_1,t')$. Let $s=2t+t'$ and $\varphi$ be an admissible trajectory, defined on $[0,s]$ and that verifies
$$
\varphi(0)=\varphi(s)=x,\qquad \varphi(t)=x_1, \qquad \varphi(t+t')=x_2.
$$
If $\tau\in[0,s]$ and $y=\varphi(\tau)$, then $x\in \A(y,s-\tau)$ and $x_1\in \Int(\A(y,s-\tau+t))$. Since $2s\geq s-\tau+t$ we obtain
$y\in \Int(\A(y,2s))$
and, considering the time-reversed system,
$$
y\in \ \ \Int(\A(y,2s))\bigcap \Int(\A^-(y,2s)).
$$
Consequently the open sets
$$
V_{\tau}=\Int(\A(\varphi(\tau),2s))\bigcap \Int(\A^-(\varphi(\tau),2s)),
$$
defined for $\tau \in [0,s]$, cover the set $K=\varphi([0,s])$. Since this last is compact we can find $0=t_0<t_1<\dots <t_k<s$ such that the
$V_{t_i}$'s, $i=0,\dots , k$ cover $K$.

Let $i\neq j$ such that $V_{t_i}\bigcap V_{t_j}\neq \emptyset$, and $z\in V_{t_i}\bigcap V_{t_j}$. 
We have $z\in \Int(\A(\varphi(t_i),2s))\bigcap \Int(\A^-(\varphi(t_j),2s))$ hence $\varphi(t_j)\in \A(\varphi(t_i),4s)$.
Therefore $\varphi(t_i)\in \A(x,4ks)$ for $i=0,\dots ,k$ and, as any point $y\in K$ belongs to one of the $V_{t_i}$'s, we obtain $K\subset \Int(\A(x,(4k+2)s))$.

Let $S=(4k+2)s$ and $t\geq S$. We can extend $\varphi$ to $\R$ by periodicity, and set $y=\varphi(-(t-S))$. Then $y\in \Int(\A(S,x))$ (because $y\in K$), and $x\in \A(y,t-S)$. Finally $x\in \Int(\A(x,t))$, and the proof is finished.

\hfill $\Box$

\vskip 0.2 cm

\noindent\textit{Proof of Theorem \ref{exactime}}

The manifold $N$ being compact, and according to Theorem \ref{finitetime}, there exists $T>0$ such that the system is controllable in time less than or equal to $T$.

Let $x$ and $S$ as in Lemma \ref{Thelemma}, and $t\geq S+2T$.

Let $y$ and $z$ two points of $N$. There exists $t_1$ (resp. $t_2$), with $0\leq t_1\leq T$ (resp. $0\leq t_2\leq T$), such that $x\in \A(y,t_1)$ and $z\in \A(x,t_2)$. Let $\tau=t-t_1-t_2$. As $\tau\geq S$ we know by Lemma \ref{Thelemma} that $x\in \A(x,\tau)$. Therefore $z\in \A(y,t)$, and the points $y$ and $z$ being arbitrary, the system is controllable at exact time $t$.

In conclusion the time $t_0$ can be chosen equal to $S+2T$.

For the converse assume the system to be Lie-determined. Then
$$
\forall x\in N \qquad \rank \Li(x)=n=\dim(N).
$$
If the rank of $\Lo$ is nowhere equal to $n$, then it is constant, equal to $n-1$, and the distribution spanned by $\Lo$ is integrable (notice that it is by construction involutive). As well as $\Li$ the family $\Lo$ is therefore Lie-determined. The orbit $\mathcal{O}_x$ of this distribution through a point $x$ is the so-called zero-time orbit of $x$ and is a $(n-1)$-dimensional submanifold of $N$. But in that case the set $\A(x,t)$ is for all $t\geq 0$ included in the $t$-translate of $\mathcal{O}_x$ (see \cite{Jurdjevic97}, Theorem 3, page 71)(the state space being compact, the requirement that at least one vector field of the family is complete is satisfied). This implies that $\Int(\A(x,t))=\emptyset$ for all $t>0$ and all $x\in N$.

\hfill $\Box$

Another sufficient condition for exact time controllability is that at least one vector field of the system vanishes at some point.

\newtheorem{singularity}[finitetime]{Theorem}
\begin{singularity}\label{singularity}
The manifold $N$ is assumed to be compact and the system $(\Sigma)$ to be controllable. If there exists one point $x$ and one control value $u\in\U$ such that $f(x,u)=0$ then there exists $t_0>0$ such that $(\Sigma)$ is controllable at exact time $t$ for all $t\geq t_0$.
\end{singularity}

\noindent\textit{Proof of Theorem \ref{singularity}}

Since $f(x,u)=0$, an admissible trajectory can stay at $x$ for some time. Therefore
$$
\forall t\geq 0 \qquad \A(x,t)=\A(x,\leq t) \quad\mbox{and}\quad \A^-(x,t)=\A^-(x,\leq t)
$$
Let $T>0$ such that $(\Sigma)$ is controllable in time less or equal to $T$. Then
$$
\forall \tau \geq T \quad \A(x,\tau)=\A(x,\leq \tau)=N \quad\mbox{and}\quad \A^-(x,\tau)=\A^-(x,\leq \tau)=N.
$$
It follows that $(\Sigma)$ is controllable at exact time $t$ for all $t\geq 2T$.

\hfill $\Box$

\vskip 0.2cm

\noindent{\bf Remark.} In the proofs of both Theorems \ref{exactime} and \ref{singularity} the compactness assumption is used only to ensure that controllability implies finite time controllability. Therefore their conclusions apply to finite time controllable systems on noncompact manifolds, and we can state:
\newtheorem{noncompactbis}[finitetime]{Theorem}
\begin{noncompactbis}\label{noncompactbis}
Let $(\Sigma)$ be a finite time controllable system. It is exact time controllable as soon as one of the conditions

(i) There exists a point $x$ where $\rank(\Lo)(x)=\dim N$;

(ii) There exist a point $x$ and a control value $u\in\U$ such that $f(x,u)=0$;

\noindent is satisfied.

If the system is Lie-determined, then Condition (i) is also necessary.
\end{noncompactbis}


\section{Applications}\label{Applications}

Before looking at some examples, notice that Theorem \ref{finitetime} (resp. Theorems \ref{exactime} and \ref{singularity}) does not say that a system satisfying its hypothesis is controllable in time less than or equal to $t$ (resp. at exact time $t$) {\it for any $t>0$}. 

It is shown in \cite{JS272} that a right invariant system on a compact and connected semisimple Lie group is exact time controllable as soon as it satisfies the rank condition. However the authors exhibit such a controllable system on the group $SO_3$ that cannot be controlled in arbitrary small time.

\subsection{Even-dimensional spheres}

It is well known that on even-dimensional spheres, any vector field vanishes at least at one point. Consequently a controllable system on an even-dimensional sphere is firstly finite time controllable by Theorem \ref{finitetime}, and then exact time controllable according to Theorem \ref{singularity}.

\subsection{Linear systems on compact Lie groups}
A linear system on a connected Lie group $G$ is a control affine one
$$
(L) \qquad\qquad\qquad \dot{x}=\X(x)+\sum_{j=1}^m u_jY_j(x)
$$
where the $Y_j$'s are right invariant vector fields and $\X$ is linear. This means that the flow of $\X$ is a one parameter group of automorphisms
(see for instance \cite{JouanAX08}).

Whenever the group $G$ is compact such a system is controllable if and only if it satisfies the rank condition (see \cite{AS00} or \cite{CM05}). As a linear vector field vanishes at the identity $e$ (because $e$ is a fixed point of any automorphism) Theorem \ref{singularity} applies. Since $\X(e)=0$ we have also $\rank \Lo(e)=\rank \Li(e)=\dim N$ as soon as the rank condition is satisfied and Theorem \ref{exactime} applies as well. We can therefore state:

\newtheorem{linear}[finitetime]{Theorem}
\begin{linear} \label{linear}
A linear system on a compact and connected Lie group is controllable if and only if it satisfies the rank condition. It is in that case exact time controllable.
\end{linear}

A vector field on a connected Lie group is said to be \textit{affine} if it is equal to the sum of a linear vector field and a right-invariant one (see \cite{JouanAX08}). The system obtained by replacing $\X$ by an affine vector field $F$ in $(L)$ may again be called linear, but Theorem \ref{linear} is no longer true in that case. Indeed the drift vector field $F$ may be merely right invariant (that is the linear part is zero). However controllable but not exact time controllable right-invariant systems on compact Lie groups are known (for instance on Torus).

In fact linear vector fields on Lie groups are analytic, as well as left or right invariant ones, hence linear systems are Lie-determined. Their exact time controllability depends upon the rank of $\Lo$.

\subsection{Analytic systems}

Analytic systems have been known to be Lie determined for a long time (see \cite{JS172} and \cite{Jurdjevic97}), so that following Theorem \ref{exactime} a controllable analytic system on a compact and connected manifold is exact time controllable if and only if the rank of $\Lo$ is full. In the famous paper \cite{JS172} the authors proved that an analytic system cannot be controllable without having the property that the rank of $\Lo$ is everywhere full in the two following cases
\begin{enumerate}
\item the covering space of the manifold $N$ is compact;
\item the fundamental group of $N$ has no elements of infinite order.
\end{enumerate}
From this we can for instance deduce that an analytic system on a compact and simply connected manifold (for example a sphere $S^n$ with $n\geq 2$) is exact time controllable if and only if it is controllable.



\begin{thebibliography}{99}
\bibitem[AS00]{AS00} V.~Ayala and L.~San Martin \textit{Controllability properties of a class of control systems on Lie groups}, Nonlinear control in the year 2000, Vol. 1 (Paris),  83--92, LN in Control and Inform. Sci., 258, Springer, 2001.
\bibitem[CM05]{CM05} F.~Cardetti and D.~Mittenhuber \textit{Local controllability for linear control systems on Lie groups}, Journal of Dynamical and Control Systems, Vol. 11, No. 3, July 2005, 353-373.
\bibitem[Grasse84]{Grasse84} K.A.~Grasse \textit{On accessibility and normal accessibility; the openness of controllability in the fine $\cc^0$ topology}, J. Differential equations, 53 (1984), 387-414.
\bibitem[GS90]{GS90} K.A.~Grasse and  H.J.~Sussmann \textit{Global Controllability by nice controls} in Nonlinear Controllability and Optimal control, Sussmann editor, Dekker, New-York, 1990, 33-79.
\bibitem[JouanAX08]{JouanAX08} Ph.~Jouan \textit{Equivalence of Control Systems with Linear Systems on
  Lie Groups and Homogeneous Spaces} arXiv 0812.0058, November 2008.
\bibitem[JS172]{JS172} V.~Jurdjevic and H.J.~Sussmann \textit{Controllability of Nonlinear Systems}, Journal of Differential Equations 12, 95-116 (1972).
\bibitem[JS272]{JS272} V.~Jurdjevic and H.J.~Sussmann \textit{Control Systems on Lie Groups}, Journal of Differential Equations 12, 313-329 (1972).
\bibitem[Jurdjevic97]{Jurdjevic97} V.~Jurdjevic \textit{Geometric control theory}, Cambridge university press, 1997.
\bibitem[Sontag98]{Sontag98} E.D.~Sontag \textit{Mathematical Control Theory}, second edition,\\ Springer, New-York, 1998.
\bibitem[Sussmann73]{Sussmann73} H.J.~Sussmann \textit{Orbits of families of vector fields and integrability of distributions}, Trans. Amer. Math. Soc. 180, 1973, 171-188.
\bibitem[Sussmann76]{Sussmann76} H.J.~Sussmann \textit{Some properties of vector field systems that are not altered by small perturbations}, J. Differential equations, 20 (1976), 292-315.
\end{thebibliography}
\end{document}